\documentclass[12pt]{amsart}
\usepackage{amssymb,amsmath,amsthm,latexsym}
\usepackage{mathrsfs}
\usepackage{mdwlist}
\usepackage{a4wide}

\newtheorem*{lemma1}{Lemma 1}
\newtheorem*{lemma2}{Lemma 2}

\theoremstyle{definition}

\theoremstyle{remark}

\newcounter{smallromans}

\newenvironment{romanenumerate}
{\begin{list}{{\normalfont\textrm{(\roman{smallromans})}}}%
    {\usecounter{smallromans}\setlength{\itemindent}{0cm}%
      \setlength{\leftmargin}{5.5ex}\setlength{\labelwidth}{5.5ex}%
      \setlength{\topsep}{0.75\parsep}\setlength{\partopsep}{0ex}%
      \setlength{\itemsep}{0ex}}}%
  {\end{list}}

\newcounter{smallalphs}

  {\end{list}}

\author{Tomasz Kania}
\address{Department of Mathematics and Statistics, Fylde College,
  Lancaster University, Lancaster LA1 4YF, United Kingdom \newline
Current address: Institute of Mathematics, Polish Academy of Sciences, \'{S}niadeckich~8, 00-956 Warszawa, Poland}
\email{tomasz.marcin.kania@gmail.com}
\title{A short proof of the fact that the matrix trace is the expectation of the numerical values}
\begin{document}
\maketitle
\begin{abstract}Using the fact that the normalised matrix trace is the unique linear functional $f$ on the algebra of $n\times n$ matrices which satisfies $f(I)=1$ and $f(AB)=f(BA)$ for all $n\times n$ matrices $A$ and $B$, we derive a well-known formula expressing the normalised trace of a complex matrix $A$ as the expectation of the numerical values of $A$; that is the function $\langle Ax,x\rangle$, where $x$ ranges the unit sphere of $\mathbb{C}^n$.\end{abstract}
Let $A=[a_{ij}]$ be an $n\times n$ complex matrix. The aim of this note is to give an easy proof of the fact that the normalised trace of $A$, $\mbox{tr}\,A = \tfrac{1}{n}(a_{11} + a_{22} + \ldots + a_{nn})$, can be thought of as the expectation of the numerical values of $A$; that is the function $x\mapsto \langle Ax,x\rangle$ defined on the Euclidean unit sphere in $\mathbb{C}^n$, endowed with the normalised Lebesgue surface measure $\mu$. More precisely,
\begin{equation}\label{main} \mbox{tr}\,A = \int\limits_{\|x\|=1} \langle Ax, x\rangle \,\mu(\mbox{d}x).\end{equation}
The above formula is a particular version of a more general identity for symmetric 2-tensors on Riemannian manifolds (consult \textit{e.g.} 
\cite{riem}; see also \cite{chow} for the proof\footnote{Bennett Chow blames one of his students for a neat proof he posted at \cite{chow}.}). We offer here an elementary proof of Equation~\eqref{main} relying on two folklore facts from linear algebra. 
\begin{lemma1}\label{unique}The matrix trace is unique in the sense that it is the unique linear functional $f\colon M_n(\mathbb{C})\to \mathbb{C}$ satisfying the following properties:
\begin{romanenumerate}
\item \label{i} $f(I)= 1$, where $I$ denotes the $n\times n$ identity matrix,
\item \label{ii} $f(AB) = f(BA)$ for all $A,B\in M_n(\mathbb{C})$.
\end{romanenumerate}
\end{lemma1}\noindent
It is evident that the standard normalised trace $\mbox{tr}$ on $M_n(\mathbb{C})$ satisfies conditions \eqref{i}-\eqref{ii}, so in order to prove the above lemma, it is enough to show that a functional $f$ enjoying \eqref{i}-\eqref{ii} agrees with $\mbox{tr}$ on the standard matrix units $e_{ij} = [\delta_{i,j}]\;(1\leqslant i,j\leqslant n)$; that is, $f(e_{ij}) = \tfrac{1}{n}$ whenever $i=j$ and $f(e_{ij}) = 0$ otherwise. We leave this as an exercise for the reader.

We shall require also the following easy and well-known fact. (See also \cite[Lemma 3.2.21]{ped}.)
\begin{lemma2}Every complex $n\times n$ matrix is a linear combination of unitary matrices.\end{lemma2}
\begin{proof}Every matrix $A\in M_n(\mathbb{C})$ can be written as a linear combination of two self-adjoint matrices, so without loss of generality it is enough to show that each self-adjoint matrix $A$ with $\|A\|\leqslant 1$ can be written as a linear combination of unitaries. To this end, set $U = A - i(I-A^2)^{\frac{1}{2}}$ and note that $U$ is unitary. Clearly $A=\frac{1}{2}U + \frac{1}{2}U^*$.\end{proof}

We are now in a position to prove Equation \eqref{main}.
\begin{proof}Let $f(A)$ denote the right hand side of Equation \eqref{main}. It is enough to verify that $f$ meets conditions \eqref{i}-\eqref{ii} of Lemma~1. Evidently, $f$ is linear, and $f(I)=1$ because $\mu$ is a probability measure. It remains to show that \eqref{ii} holds.

Let $A,B\in M_n(\mathbb{C})$ and let us write $B$ as a linear combination of some unitary matrices $U_1, \ldots, U_m$, that is $B=\sum_{k=1}^m a_k U_k$ for some scalars $a_1, \ldots, a_m$. We may assume additionally that each matrix $U_k \; (k\leqslant n)$ has determinant 1, as we can always write $B=\sum_{k=1}^m (a_k\det U_k) \frac{U_k}{\det U_k}$. We have $f(AB)=f(BA)$ as soon as $f(AU_k)=f(U_kA)$ for all $k\leqslant m$, so that without lost of generality we may suppose that $B$ is unitary and $\det B =1$. Making the substitution $x=B^*z$ and taking into account that the determinant of $B$ is equal to 1 (hence also the Jacobian of $B$, regarded as a map from the real $(2n-1)$-sphere to itself, is equal to 1), we arrive at the conclusion that
\[\begin{array}{lcl}f(AB) &=& \int\limits_{\|x\|=1} \langle ABx, x\rangle \,\mu(\mbox{d}x)\\
 &=&  \int\limits_{\|x\|=1} \langle Az, B^*z \rangle \,\mu(\mbox{d}z) \\
 &=& \int\limits_{ \|x\|=1} \langle BAz, z \rangle \,\mu(\mbox{d}z) \\
&=& f(BA),\end{array}\]
which completes the proof.\end{proof}

\bibliographystyle{amsplain}

\end{document}